\documentclass[a4paper,oneside,11pt]{article}
\usepackage[a4paper]{geometry}
\usepackage{aeguill}                 
\usepackage{graphicx}
\usepackage{caption}
\usepackage{amsmath,amsfonts,amssymb,amsthm}
\usepackage{mathabx}
\usepackage{pstricks,comment} 
\usepackage{pst-plot}
\usepackage{pstricks-add}
\usepackage{multido}
\usepackage{url}
\usepackage{epstopdf,enumerate}
\usepackage{srcltx}
\usepackage[utf8]{inputenc} 
\usepackage[T1]{fontenc} 
\usepackage[english]{babel} 
\usepackage{hyperref}
\usepackage[all]{xy} 
\usepackage{titlesec}
\usepackage{mathabx}
\usepackage{tikz}
\usepackage{fancyhdr}
\usepackage{mathrsfs}
\usepackage[us,12hr]{datetime}
\fancypagestyle{plain}{
\fancyhf{}
\rfoot{ {\ddmmyyyydate\today}}
\lfoot{\thepage}
}
\renewcommand{\theenumi}{\alph{enumi}}

    \newtheoremstyle{TheoremNum}
        {\topsep}{\topsep}              
        {\itshape}                      
        {}                              
        {\bfseries}                     
        {.}                             
        { }                             
        {\thmname{#1}\thmnote{ \bfseries #3}}

{\theoremstyle {definition} \newtheorem {defi} {Definition}[section]}
{\theoremstyle {plain}  \newtheorem {theo} [defi] {Theorem}}
{\theoremstyle {plain}  \newtheorem {coro} [defi] {Corollary}}
{\theoremstyle {plain} \newtheorem {prop} [defi] {Proposition}}
{\theoremstyle {plain} \newtheorem {lem}[defi] {Lemma}}
{\theoremstyle {plain} }
{\theoremstyle {plain} }
{\theoremstyle{TheoremNum} }
{\theoremstyle{TheoremNum} }
{\theoremstyle{TheoremNum} }

\newcommand{\Aut}{\mathrm{Aut}}
\newcommand{\Out}{\mathrm{Out}}

\newcommand{\Stab}{\mathrm{Stab}}

\newcommand{\Fix}{\mathrm{Fix}}
\newcommand{\IA}{\mathrm{IA}}

\newcommand{\id}{\mathrm{id}}
\newcommand{\ad}{\mathrm{ad}}
\newcommand{\FF}{\mathrm{FF}}
\newcommand{\ZZ}{\mathbb{Z}}

\newcommand{\NN}{\mathbb{N}}
\newcommand{\RR}{\mathbb{R}}

\newcommand{\dem}{\noindent{\bf Proof. }}

\title{Roots of outer automorphisms of free groups and centralizers of abelian subgroups of $\Out(F_N)$}
\author{Yassine Guerch}
\date{\today}
\begin{document}
\maketitle
\renewcommand*\labelenumi{(\theenumi)}

\begin{abstract}
Let $N \geq 2$ and let $\Out(F_N)$ be the outer automorphism group of a nonabelian free group of rank $N$. Let $\IA_N(\ZZ/3\ZZ)$ be the finite index subgroup of $\Out(F_N)$ which is the kernel of the natural action of $\Out(F_N)$ on $H_1(F_N,\ZZ/3\ZZ)$. We show that $\IA_N(\ZZ/3\ZZ)$ is an \emph{R-group}, that is, for every $\phi,\psi \in \IA_N(\ZZ/3\ZZ)$, if there exists $k \in \NN^*$ such that $\phi^k=\psi^k$, then $\phi=\psi$. This answers a question of Handel and Mosher. We then use the fact that $\IA_N(\ZZ/3\ZZ)$ is an $R$-group in order to prove that the normalizer in $\IA_N(\ZZ/3\ZZ)$ of every abelian subgroup of $\IA_N(\ZZ/3\ZZ)$ is equal to its centralizer. We finally give an alternative proof of a result, due to Feighn and Handel, that the centralizer of an element of $\Out(F_N)$ which has only finitely many periodic orbits of conjugacy classes of maximal cyclic subgroups of $F_N$ is virtually abelian.
\end{abstract}

\section{Introduction}\label{Section Introduction}

Let $N \geq 2$ and let $\Out(F_N)$ be the outer automorphism group of a nonabelian free group $F_N$. When studying dynamical or algebraic properties of $\Out(F_N)$, one often needs to pass to a finite index subgroup in order to avoid periodic behaviours. A finite index subgroup of particular importance for such considerations is the subgroup $\IA_N(\ZZ/3\ZZ)$, which is the kernel of the natural action of $\Out(F_N)$ on the first homology group of $F_N$ with coefficients in $\ZZ/3\ZZ$. The group $\IA_N(\ZZ/3\ZZ)$ is an analogue of the congruence subgroups of arithmetic lattices, which are finite index subgroups of the considered lattices having the advantage of being torsion free. The group $\IA_N(\ZZ/3\ZZ)$ similarly satisfies numerous aperiodic properties. For instance, it is torsion free~\cite[Corollary~5.7.6]{BesFeiHan00}, every virtually abelian subgroup of $\IA_N(\ZZ/3\ZZ)$ is abelian~\cite{HandelMosher20abelian}, and any conjugacy class of a free factor of $F_N$ which has a periodic orbit by a subgroup of $\IA_N(\ZZ/3\ZZ)$ is in fact fixed~\cite[Theorem~II.3.1]{HandelMosher20}. We notice that we could have replaced $\ZZ/3\ZZ$ in the definition of $\IA_N(\ZZ/3\ZZ)$ by any other $\ZZ/m\ZZ$ with $m \geq 3$ without changing the aperiodic properties of the group. On the contrary, the group $\IA_N(\ZZ/2\ZZ)$ is not torsion free as it contains the involution sending every element of some basis of $F_N$ to its inverse.

This is why the group $\IA_N(\ZZ/3\ZZ)$ is predominant when studying rigidity properties of $\Out(F_N)$. It is for instance of central importance in the proof of the Tits alternative due to Bestvina, Feighn and Handel~\cite{BesFeiHan00,BesFeiHan05}. The fact that periodic conjugacy classes of free factors and elements of $F_N$ are fixed by subgroups of $\IA_N(\ZZ/3\ZZ)$ is a key step in the proof of decomposition theorems due to Handel and Mosher~\cite{HandelMosher20} and in their alternative regarding the second bounded cohomology
of subgroups of $\Out(F_N)$~\cite{HandelMosher2017}. It also intervenes in the proof of the rigidity of the abstract commensurator of $\Out(F_N)$ due to Farb and Handel~\cite{farb2007commensurations} for $N \geq 4$ and Horbez and Wade~\cite{HorbezWade20} for $N \geq 3$. Finally, it also appears in the recent proof of the measure equivalence rigidity of $\Out(F_N)$ due to Guirardel and Horbez~\cite{Guirardelhorbez}.

In this article, we prove another aperiodic property of $\IA_N(\ZZ/3\ZZ)$, which is the content of the following theorem.

\begin{theo}\label{Theo intro}
Let $N \geq 2$ and let $\phi,\psi \in \IA_N(\ZZ/3\ZZ)$. Suppose that there exists $k \in \NN^*$ such that $\phi^k=\psi^k$. Then $\phi=\psi$.
\end{theo}

This answers a question of Handel and Mosher~\cite[Question~1.4]{HandelMosher20abelian}. Following Kontorovi\v{c}~\cite{Kontorovic1946} (see also~\cite{Baumslag1960}), Theorem~\ref{Theo intro} implies that $\IA_N(\ZZ/3\ZZ)$ is an \emph{$R$-group}. This implies in particular (see Corollary~\ref{Coro R group}) that the normalizer in $\IA_N(\ZZ/3\ZZ)$ of any cyclic subgroup of $\IA_N(\ZZ/3\ZZ)$ is equal to its centralizer. 

Another consequence of Theorem~\ref{Theo intro} is that one can enlarge the recognition theorem due to Feighn and Handel~\cite[Theorem~5.3]{FeiHan06} to a larger class of elements of $\Out(F_N)$. Indeed, the recognition theorem only applies to \emph{forward rotationless} elements. By~\cite[Lemma~4.42]{FeiHan06}, there exists $K_N \in \NN^*$ such that, for every $\phi \in \Out(F_N)$, the element $\phi^{K_N}$ is forward rotationless. Thus, we have the following corollary.

\begin{coro}\label{Coro Recognition theorem}
Let $N \geq 2$  and let $K_N \in \NN^*$ be the constant given by~\cite[Lemma~4.42]{FeiHan06}. Let $\phi,\psi \in \IA_N(\ZZ/3\ZZ)$. Then $\phi=\psi$ if and only if $\phi^{K_N}$ and $\psi^{K_N}$ satisfy the hypotheses of the recognition theorem~\cite[Theorem~5.3]{FeiHan06}.
\end{coro}

Recall that, in an $R$-group, the normalizer of every cyclic subgroup is equal to its centralizer. We may ask whether this property extends to every \emph{abelian} subgroup of $\IA_N(\ZZ/3\ZZ)$. This is in fact the case as shown by the following theorem.

\begin{theo}[see~Theorem~\ref{Theo Normalizer equals centralizer}]\label{Theo normalizer equal centralizer intro}
Let $N \geq 2$. For every abelian subgroup $H \subseteq \IA_N(\ZZ/3\ZZ)$, the normalizer of $H$ in $\IA_N(\ZZ/3\ZZ)$ is equal to its centralizer.
\end{theo} 

Note that Theorem~\ref{Theo normalizer equal centralizer intro} is no longer true if we remove the assumption of being abelian. We thus deduce this new algebraic characterization of abelian subgroups of $\IA_N(\ZZ/3\ZZ)$.

\begin{coro}\label{Characterization abelian subgroups}
Let $H$ be a subgroup of $\IA_N(\ZZ/3\ZZ)$. Then $H$ is abelian if and only if its normalizer in $\IA_N(\ZZ/3\ZZ)$ equals its centralizer.
\end{coro}

Theorem~\ref{Theo normalizer equal centralizer intro} extends a similar result due to Handel and Mosher for abelian \emph{UPG} subgroups of $\IA_N(\ZZ/3\ZZ)$~\cite[Proposition~1.3]{HandelMosher20abelian}. It also enlarges the known informations regarding both abelian subgroups of $\Out(F_N)$ (see for instance~the work of Feighn and Handel~\cite{FeiHan09}) and centralizers of subgroups of $\Out(F_N)$, see for instance the work of Rodenhausen and Wade~\cite{RodenhausenWade2015}, Algom-Kfir and Pfaff~\cite{AlgomKfirPfaff2017} and Andrew and Martino~\cite{AndrewMartino2022}. 

As we will see in Corollary~\ref{Coro obstruction}, Theorem~\ref{Theo normalizer equal centralizer intro} gives obstructions to the existence of subgroups of $\IA_N(\ZZ/3\ZZ)$ isomorphic to some semidirect products of groups. Indeed, if $H$ is a subgroup of $\IA_N(\ZZ/3\ZZ)$ isomorphic to $A \rtimes B$ where $A$ is abelian, then Theorem~\ref{Theo normalizer equal centralizer intro} implies that $H$ is in fact isomorphic to the direct product of $A$ and $B$.

Finally, our techniques also enable us to understand centralizers of some elements of $\Out(F_N)$. The following theorem also follows from the work of Feighn and Handel~\cite{FeiHan09} regarding abelian subgroups of $\Out(F_N)$ (see also the work of Mutanguha~\cite{Mutanguha22} regarding centralizers of \emph{atoroidal} elements of $\Out(F_N)$).

\begin{theo}[Theorem~\ref{Prop atoroidal}]\cite{FeiHan09}\label{Theo intro atoroidal}
Let $N \geq 2$ and let $\phi \in \Out(F_N)$ be an outer automorphism which has only finitely many periodic orbits of conjugacy classes of maximal cyclic subgroups of $F_N$. The centralizer of $\phi$ in $\Out(F_n)$ is virtually abelian.
\end{theo}

We note that the work of Feighn and Handel uses the technology of train tracks, while the proof of Mutanguha uses some actions on limit trees. Our proof relies on isometric actions of $\Out(F_N)$ on some Gromov-hyperbolic spaces, so that the techniques used in this paper significantly differs from the other two proofs.

We now sketch the proof of Theorem~\ref{Theo intro}, the proofs of Theorems~\ref{Theo normalizer equal centralizer intro} and~\ref{Theo intro atoroidal} following essentially the same lines. Let $\phi, \psi \in \IA_N(\ZZ/3\ZZ)$ and let $k \in \NN^*$ be such that $\phi^k=\psi^k$. The idea of the proof is to show that the group $H=\langle \phi,\psi \rangle$ is abelian. Since $\IA_N(\ZZ/3\ZZ)$ is torsion free and since $\phi^k=\psi^k$, this will imply that $H$ is cyclic and this will conclude the proof. 

The proof is by induction on $N$, the case $N=1$ being immediate. Consider a maximal proper $H$-invariant \emph{free factor system $\mathcal{F}$ of $F_N$} (see Section~\ref{Subsection malnormal}).  If $\mathcal{F}$ is \emph{sporadic}, then one can canonically associate to $\mathcal{F}$ a Bass-Serre tree $S$ of $F_N$ whose conjugacy classes of vertex stabilizers are contained in $\mathcal{F}$. The action of $F_N$ on $S$ naturally extends to an action of the preimage $\widetilde{H}$ of $H$ in $\Aut(F_N)$. A thorough investigation of the action of $\widetilde{H}$ then shows that $H$ is necessarily abelian. This investigation relies on the induction hypothesis in order to deal with the action of $\widetilde{H}$ on vertex stabilizers of $S$ in $F_N$.

Suppose now that $\mathcal{F}$ is \emph{nonsporadic}. Then $H$ acts by isometries on a Gromov hyperbolic space called the \emph{free factor graph of $F_N$ relative to $\mathcal{F}$} and denoted by $\FF(F_N,\mathcal{F})$. Maximality of $\mathcal{F}$ implies that $H$ contains a loxodromic element by results of Handel and Mosher~\cite[Theorem~A]{HandelMosher20} and Guirardel and Horbez~\cite[Theorem~7.1]{Guirardelhorbez19}. Moreover, since $H$ is abelian, it acts with a finite orbit on the Gromov boundary $\partial_{\infty}\FF(F_N,\mathcal{F})$ of $\FF(F_N,\mathcal{F})$. A thorough analysis of stabilizers of points in $\partial_{\infty}\FF(F_N,\mathcal{F})$ based on their description by Guirardel and Horbez~\cite{Guirardelhorbez19} and Horbez and Wade~\cite{HorbezWade20} then concludes the proof. The analysis is essentially done in Proposition~\ref{Prop nonsporadic case}.

We now present the structure of the article. In Section~\ref{Section preliminaries}, we recall basic definitions associated with $\Out(F_N)$ such as free factor systems, splittings and properties of the action of $\Out(F_N)$ on the relative free factor graph. Sections~\ref{Section Sporadic} and~\ref{Section Proof} are devoted to the proof of Theorem~\ref{Theo intro}. Section~\ref{Section Sporadic} focuses on the sporadic case and we finish the proof of Theorem~\ref{Theo intro} in Section~\ref{Section Proof}. Finally, in Section~\ref{Section abelian}, we prove Theorem~\ref{Theo normalizer equal centralizer intro} regarding the normalizer of abelian subgroups of $\IA_N(\ZZ/3\ZZ)$ and Theorem~\ref{Theo intro atoroidal} regarding the centralizer of some elements of $\Out(F_N)$.

\medskip

{\small {\bf Acknowledgments. } I warmly thank Damien Gaboriau, Camille Horbez and Frédéric Paulin for their precious advices as well as their very helpful remarks regarding earlier versions of this work. I am also grateful to Jean Pierre Mutanguha for his numerous very useful comments and for pointing out an improvement of an earlier version of Theorem~\ref{Theo intro atoroidal}.

The author was supported by the LABEX MILYON of Université de Lyon.}

\section{Preliminaries}\label{Section preliminaries}

\subsection{Free factor systems of $F_{N}$}\label{Subsection malnormal}

Let $N \geq 2$ and let $F_N$ be a nonabelian free group of rank $N$. In this section, we present the definition of free factor systems of $F_N$, which are specific finite sets of conjugacy classes of subgroups of $F_N$. First, we give the definition of some larger collections of conjugacy classes of subgroups of $F_N$ called \emph{subgroup systems}.

\begin{defi}[Subgroup system]
A \emph{subgroup system} is a finite (possibly trivial) set $\mathcal{A}$ of conjugacy classes of finitely generated subgroups of $F_N$.
\end{defi}

There exists a natural partial order on the set of subgroup systems of $F_N$, where $\mathcal{A}_1 \leq \mathcal{A}_2$ if for every subgroup $A$ of $F_N$ such that $[A] \in \mathcal{A}_1$, there exists $[B] \in \mathcal{A}_2$ such that $A \subseteq B$. The subgroup system $\mathcal{A}_2$ is then an \emph{extension of $\mathcal{A}_1$}.

The group $\Out(F_N)$ has a natural action on the set of subgroup systems, and this action preserves the partial order defined above. Given a subgroup system $\mathcal{A}$, we denote by $\Out(F_N,\mathcal{A})$ the subgroup of $\Out(F_N)$ preserving $\mathcal{A}$.

Let $\mathcal{A}$ be a subgroup system of $F_N$. An element $g$ of $F_N$ is \emph{$\mathcal{A}$-peripheral} if there exists $[A] \in \mathcal{A}$ such that $g \in A$. Otherwise, we say that $g$ is $\mathcal{A}$-nonperipheral. A subgroup of $F_N$ is $\mathcal{A}$-nonperipheral if it contains an $\mathcal{A}$-nonperipheral element, and is $\mathcal{A}$-peripheral otherwise.

We now present an important family of subgroup systems, the \emph{free factor systems}. A \emph{free factor system} is a subgroup system $\mathcal{F}=\{[A_1],\ldots,[A_k]\}$ of $F_N$ such that there exists a subgroup $B$ of $F_N$ with $F_N=A_1 \ast \ldots \ast A_k \ast B$. 

\begin{defi}[Sporadic extension]
Let $\mathcal{F}_1 \leq \mathcal{F}_2$ be two free factor systems of $F_N$. The extension $\mathcal{F}_1 \leq \mathcal{F}_2$ is \emph{sporadic} if one of the following holds:

\medskip

\noindent{$(1)$ } there exist $[A],[B] \in \mathcal{F}_1$ such that $\mathcal{F}_2=(\mathcal{F}_1-\{[A],[B]\}) \cup \{[A \ast B]\}$;

\medskip

\noindent{$(2)$ } there exist $[A] \in \mathcal{F}_1$ and $g \in F_N$ such that $\mathcal{F}_2=(\mathcal{F}_1-\{[A]\}) \cup \{[A \ast \langle g \rangle]\}$;

\medskip

\noindent{$(3)$ } there exists $g \in F_N$ such that $\mathcal{F}_2=\mathcal{F}_1 \cup \{[\langle g \rangle]\}$.

\medskip

Otherwise, we say that the extension $\mathcal{F}_1 \leq \mathcal{F}_2$ is \emph{nonsporadic}.
\end{defi}

A free factor system $\mathcal{F}$ of $F_N$ is \emph{sporadic} if the extension $\mathcal{F} \leq \{[F_N]\}$ is sporadic. Note that, for a sporadic free factor system $\mathcal{F}$, either $\mathcal{F}=\{[A]\}$ or $\mathcal{F}=\{[A],[B]\}$ for some subgroups $A,B \subseteq F_N$.

Let $\phi \in \Out(F_N)$ and let $\mathcal{F}$ be a free factor system of $F_N$. Suppose that $\phi$ fixes every element of $\mathcal{F}$. Then, for every $[A] \in \mathcal{F}$, by malnormality of $A$, the element $\phi$ induces an element $\phi|_A \in \Out(A)$.

\subsection{Splittings of $F_N$ and twists automorphisms}\label{Section splittings}

A \emph{splitting of $F_N$} is an $F_N$-equivariant homeomorphism class $\mathcal{S}$ of a minimal, simplicial action of $F_N$ on a simplicial tree $S$. Here, minimal means that $S$ does not contain a proper $F_N$-invariant subtree. If $\mathcal{F}$ is a free factor system of $F_N$, an \emph{$(F_N,\mathcal{F})$-splitting} $\mathcal{S}$ is a splitting of $F_N$ such that, for every $[A] \in \mathcal{F}$, the group $A$ is elliptic in $\mathcal{S}$. If $v$ is a vertex of $S$, we denote by $G_v$ its stabilizer. Let $V$ be a set of representatives of the $F_N$-orbits of vertices in $S$. 

A splitting $\mathcal{S}$ is \emph{free} if edge stabilizers of $S$ are trivial. Note that, in this case, the set $\{[G_v]\}_{v \in V}$ is a free factor system of $F_N$. Given a free factor system $\mathcal{F}$ of $F_N$, a \emph{Grushko $(F_N,\mathcal{F})$-free splitting} is an $(F_N,\mathcal{F})$-free splitting such that, for every $v \in V$ with nontrivial stabilizer, we have $[G_v] \in \mathcal{F}$.

The group $\Aut(F_N)$ acts on the right on the set of splittings of $F_N$ by precomposition of the action. This action passes to the quotient to give an action of $\Out(F_N)$. Note that this action preserves the set of free splittings.

Let $\mathcal{F}$ be a free factor system of $F_N$. Given two $(F_N,\mathcal{F})$-splittings $\mathcal{S}$ and $\mathcal{S}'$, we say that $\mathcal{S}$ is a \emph{refinement} of $\mathcal{S}'$ if there exist $S \in \mathcal{S}$ and $S' \in \mathcal{S}'$ such that $S'$ is obtained from $S$ by collapsing some orbits of edges. Two $(F_N,\mathcal{F})$-splittings are \emph{compatible} if they have a common refinement which is an $(F_N,\mathcal{F})$-splitting.

Let $\mathcal{S}$ be a splitting of $F_N$ and let $S$ be a representative of $\mathcal{S}$. The subgroup $\Stab(\mathcal{S})$ of elements of $\Out(F_N)$ which fix $\mathcal{S}$ has a natural description by a result of Levitt~\cite[Proposition~4.2]{levitt2005}. Indeed, there exists a natural homomorphism from $\Stab(\mathcal{S})$ to the group of graph automorphisms $\Aut_{gr}(\overline{F_N \backslash S})$ of the underlying graph of $F_N \backslash S$. We denote the kernel of this homomorphism by $K(\mathcal{S})$. Moreover, there exists a natural homomorphism $K(\mathcal{S}) \to \prod_{v \in V}\Out(G_v)$ given by the action on vertex stabilizers. The kernel of this homomorphism, denoted by $T(\mathcal{S})$, is the \emph{group of twists of $\mathcal{S}$}.

\begin{lem}\cite[Lemma~5.3]{CohenLustig99}\label{Lem Dehn twists central}
Let $\mathcal{S}$ be a splitting of $F_N$ whose edge stabilizers are all nontrivial. The group $T(\mathcal{S})$ is central in $K(\mathcal{S})$.
\end{lem}

Let $\mathcal{F}$ be a sporadic free factor system of $F_N$. Then one can naturally associate to $\mathcal{F}$ a Grushko $(F_N,\mathcal{S})$-free splitting $\mathcal{S}$, which is the \emph{Bass-Serre tree of $F_N$ associated with $\mathcal{F}$}. A representative $S \in \mathcal{S}$ has exactly one orbit of edges. The splitting $\mathcal{S}$ is fixed by $\Out(F_N,\mathcal{F})$. Moreover, by~\cite[Proposition~3.1]{levitt2005}, the group of twists of $\mathcal{S}$ is isomorphic to a direct product of two free (maybe cyclic) groups.

\subsection{Properties of the subgroup $\IA_N(\ZZ/3\ZZ)$}

Let $N \geq 2$ and let $$\IA_N(\ZZ/3\ZZ)=\ker(\Out(F_N) \to \Aut(H_1(F_N,\ZZ/3\ZZ)).$$ In this section, we recall some properties of $\IA_N(\ZZ/3\ZZ)$ that will be used in the proof of Theorem~\ref{Theo intro}. Most of the following results state the aperiodicity of some orbits associated with the action of $\IA_N(\ZZ/3\ZZ)$ on some natural sets.

\begin{prop}\cite[Corollary~5.7.6]{BesFeiHan00}\label{Prop torsion free}
The group $\IA_N(\ZZ/3\ZZ)$ is torsion free.
\end{prop}

\begin{theo}\label{Theo free factor fixed}\cite[Theorem~II.3.1]{HandelMosher20} 
Let $H$ be a subgroup of $\IA_N(\ZZ/3\ZZ)$ and let $\mathcal{F}$ be an $H$-periodic free factor system. Then $\mathcal{F}$ is fixed by $H$ and every element $[A] \in \mathcal{F}$ is fixed by $H$.
\end{theo}

\begin{theo}\label{Theo conjugacy class fixed}\cite[Theorem~II.4.1]{HandelMosher20}
Let $H$ be a subgroup of $\IA_N(\ZZ/3\ZZ)$. Then every $H$-periodic conjugacy class of some element of $F_N$ is fixed by $H$. 
\end{theo}

\begin{theo}\cite[Theorem~1.1]{HandelMosher20abelian}\label{Theo virtually abelian is abelian}
Let $H$ be a virtually abelian subgroup of $\IA_N(\ZZ/3\ZZ)$. Then $H$ is abelian and finitely generated.
\end{theo}

The fact, stated in Theorem~\ref{Theo virtually abelian is abelian}, that any abelian subgroup $H$ of $\Out(F_N)$ is finitely generated also follows from the work of Bass and Lubotzky~\cite{BassLubotzky1994}.

\begin{lem}\cite[Lemma~2.6]{HorbezWade20}\label{Lem periodic splitting is fixed}
Let $H \subseteq \IA_N(\ZZ/3\ZZ)$ and let $\mathcal{S}$ be an $H$-periodic free splitting. Then $\mathcal{S}$ is fixed by $H$ and $H \subseteq K(\mathcal{S})$.
\end{lem}

\begin{lem}\label{Lem Tits}
A subgroup $H$ of $\IA_N(\ZZ/3\ZZ)$ is abelian if and only if it does not contain a nonabelian free group.
\end{lem}

\dem By the Tits alternative due to Bestvina, Feighn and Handel~\cite{BesFeiHan00,BesFeiHan05}, either $H$ is virtually solvable, or it contains a nonabelian free group. By~\cite[Theorem~1.1]{BesFeiHan04} (see also \cite[Corollary~1.3]{Alibegovic2002}), either $H$ is virtually abelian or $H$ contains a nonabelian free group. By Theorem~\ref{Theo virtually abelian is abelian}, either $H$ is abelian or it contains a nonabelian free group.

\hfill\qedsymbol

\subsection{Relative free factor graph and relative arational trees}

Let $\mathcal{F}$ be a free factor system of $F_N$. The \emph{free factor graph of $F_N$ relative to $\mathcal{F}$}, denoted by $\FF(F_N,\mathcal{F})$, is the graph whose vertices are the $(F_N,\mathcal{F})$-free splittings of $F_N$, two such splittings being adjacent if they have a common refinement which is free or if they share a common elliptic $\mathcal{F}$-nonperipheral element.
By a result of Guirardel and Horbez~\cite[Proposition~2.11]{Guirardelhorbez19} (see also the work of Handel and Mosher~\cite{HandelMosher14} for the case $\mathcal{F}=\varnothing$), the graph $\FF(F_N,\mathcal{F})$ is Gromov-hyperbolic. In the rest of the section, we describe the group of isometries of $\FF(F_N,\mathcal{F})$ as well as its Gromov boundary.

The group $\Out(F_N,\mathcal{F})$ acts naturally on $\FF(F_N,\mathcal{F})$ by isometries. An outer automorphism $\phi \in \Out(F_N,\mathcal{F})$ is \emph{fully irreducible relative to $\mathcal{F}$} if there does not exist a proper free factor system $\mathcal{F} < \mathcal{F}'$ fixed by a power of $\phi$. These elements are in fact the loxodromic elements of $\FF(F_N,\mathcal{F})$. 

\begin{theo}\cite[Theorem~A]{gupta18}\label{Theo loxo free factor}
Let $\mathcal{F}$ be a nonsporadic free factor system of $F_N$. An element $\phi \in \Out(F_N,\mathcal{F})$ is a loxodromic element of $\FF(F_N,\mathcal{F})$ if and only if $\phi$ is fully irreducible relative to $\mathcal{F}$.
\end{theo}

Theorem~\ref{Theo loxo free factor} was proved by Gupta~\cite{gupta18} and later extended by Guirardel and Horbez to outer automorphisms of free products of groups~\cite[Theorem~4.1]{Guirardelhorbez19}. 

The following theorem was proved by Handel and Mosher~\cite{HandelMosher20} in the finitely generated case and by Guirardel and Horbez~\cite{Guirardelhorbez19} in the general case.

\begin{theo}\cite[Theorem~7.1]{Guirardelhorbez19}\cite[Theorem~A]{HandelMosher20}\label{Theo fully irreducible contained}
Let $H$ be a subgroup of $\IA_N(\ZZ/3\ZZ)$ and let $\mathcal{F}$ be a maximal proper $H$-invariant free factor system. Suppose that $\mathcal{F}$ is nonsporadic. Then $H$ contains a fully irreducible outer automorphism relative to $\mathcal{F}$.
\end{theo}

We now describe the Gromov boundary of $\FF(F_N,\mathcal{F})$. An \emph{$(F_N,\mathcal{F})$-free factor} is a subgroup of $F_N$ which arises as a point stabilizer of some $(F_N,\mathcal{F})$-free splitting. An $(F_N,\mathcal{F})$-free factor is \emph{proper} if it is $\mathcal{F}$-nonperipheral and not equal to $F_N$. Note that, if $A$ is an $(F_N,\mathcal{F})$-free factor, then $\mathcal{F}$ induces a free factor system of $A$, denoted by $\mathcal{F}|_A$. 

An \emph{$(F_N,\mathcal{F})$-arational tree} is an $\RR$-tree $T$ equipped with an $F_N$-action by isometries such that $T$ is not a Grushko $(F_N,\mathcal{F})$-free splitting and such that, for every proper $(F_N,\mathcal{F})$-free factor $A$, the action of $A$ on its minimal tree is a Grushko $(A,\mathcal{F}|_A)$-free splitting.

We record the following fact, which is a consequence of the description of the Gromov boundary of $\FF(F_N,\mathcal{F})$. It is due to Hamenstädt~\cite{Hamenstadt14} for the case $\mathcal{F}=\varnothing$, and Guirardel and Horbez~\cite{Guirardelhorbez19} for the general case.

\begin{prop}\cite[Theorem~3.4]{Guirardelhorbez19}\label{prop fix arational boundary}
Let $\mathcal{F}$ be a nonsporadic free factor system of $F_N$ and let $H$ be a subgroup of $\Out(F_N,\mathcal{F})$. If $H$ has a finite orbit in $\partial_{\infty} \FF(F_N,\mathcal{F})$, then $H$ has a finite index subgroup which fixes the homothety class of an $(F_N,\mathcal{F})$-arational tree.
\end{prop}

We now describe the stabilizer in $\Out(F_N,\mathcal{F})$ of the homothety class $[T]$ of an arational $(F_N,\mathcal{F})$-tree $T$. We have a natural homomorphism $$\mathrm{SF}\colon \Stab([T]) \to \RR_+^* $$ given by the stretching factor, whose kernel is denoted by $\Stab_{isom}(T)$. The homomorphism $\mathrm{SF}$ has the following properties. 

\begin{lem}\cite[Lemma~6.2, Proposition~6.3, Corollary~6.12]{Guirardelhorbez19}\label{Lem stratching factor cyclic}
\noindent{$(1)$ } The image of $\mathrm{SF}$ is cyclic. 

\medskip

\noindent{$(2)$ } For every $\phi \in \Stab([T])$, we have $\mathrm{SF}(\phi) \neq 1$ if and only if $\phi$ is fully irreducible relative to $\mathcal{F}$.
\end{lem}  

The following proposition is extracted from \cite{HorbezWade20}, where it is attributed to Guirardel and Levitt. Given an $(F_N,\mathcal{F})$-arational tree $T$, it describes an $(F_N,\mathcal{F})$-splitting canonically associated with a subgroup of $\Stab_{isom}(T)$. 

\begin{prop}\cite[Lemmas~5.3,~5.6, Theorem~5.4]{HorbezWade20}\label{Prop existence splitting arational tree}
Let $\mathcal{F}$ be a nonsporadic free factor system and let $T$ be an $(F_N,\mathcal{F})$-arational tree. Let $H \subseteq \Stab_{isom}(T) \cap \IA_N(\ZZ/3\ZZ)$ be a subgroup. There exists an $(F_N,\mathcal{F})$-splitting $\mathcal{S}_{T,H}$ fixed by the normalizer of $H$ in $\Out(F_N,\mathcal{F})$ with the following properties.

\medskip

\noindent{$(1)$ } Every edge stabilizer is nontrivial.

\medskip

\noindent{$(2)$ } There exists a partition $V=V_0 \amalg V_1$ of the vertices of $\mathcal{S}_{T,H}$ such that 

\begin{enumerate}
\item for every $v \in V_1$, the conjugacy class $[G_v]$ is contained in $\mathcal{F}$;

\item the set $V_0$ consists in a unique $F_N$-orbit of a vertex $v_0 \in V_0$. The group $G_{v_0}$ is a nonabelian free group. Moreover, the homomorphism
$$H \cap K(\mathcal{S}_{T,H}) \to \Out(G_{v_0}) $$ has trivial image.
\end{enumerate}

\end{prop}

We will also use the following corollary, which is an easy consequence of Proposition~\ref{Prop existence splitting arational tree}~$(2)(b)$ (see~\cite[Proposition~6.10]{Guirardelhorbez19} for a proof).

\begin{coro}\cite[Proposition~6.10]{Guirardelhorbez19}\label{Coro nonabelian fixed}
Let $\mathcal{F}$ be a nonsporadic free factor system and let $T$ be an $(F_N,\mathcal{F})$-arational tree. For every $\phi \in \Stab_{isom}(T) \cap \IA_N(\ZZ/3\ZZ)$, there exist a nonabelian free group $F$ of $F_N$ and a representative $\Phi \in \phi$ such that $\Phi(F)=F$ and $\Phi|_F=\mathrm{id}_F$.
\end{coro}

\section{The sporadic case}\label{Section Sporadic}

In this section, we prove Theorem~\ref{Theo intro} when the elements $\phi,\psi \in \IA_N(\ZZ,3\ZZ)$ fix a sporadic free factor system. We first need some preliminary lemmas.

\begin{lem}\label{Lem case automorphism}
Let $\Phi,\Psi \in \Aut(F_N)$. Suppose that there exists $k \in \NN^*$ such that $\Phi^k=\Psi^k$ and that the image $\phi$ of $\Phi$ in $\Out(F_N)$ is equal to the image $\psi$ of $\Psi$ and is contained in $\IA_N(\ZZ/3\ZZ)$. Then $\Phi=\Psi$.
\end{lem}

\dem We only need to show that $\Phi$ and $\Psi$ commute. Indeed, in this case, the group $\langle \phi,\psi \rangle$ is abelian. By Proposition~\ref{Prop torsion free}, the group $\langle \phi, \psi \rangle$ is also torsion free. Thus, since the kernel of $\Aut(F_N) \to \Out(F_N)$ is a nonabelian free group, the group $\langle \Phi,\Psi \rangle$ is free abelian. Since $\Phi^k=\Psi^k$, we have in fact $\Phi=\Psi$. 

So let us prove that $\Phi$ commutes with $\Psi$. Let $K$ be the normal subgroup of $\langle \Phi,\Psi \rangle$ consisting in all its inner automorphisms and let $K_0$ be the subgroup of $F_N$ consisting in all elements $g \in F_N$ with $\ad_g \in K$. 

We first treat the case where $K_0$ is cyclic, generated by $g \in F_N$. Since $K$ is a normal subgroup of $\langle \Phi,\Psi \rangle$, the group $K_0$ is preserved by $\langle \Phi,\Psi \rangle$. By Theorem~\ref{Theo conjugacy class fixed}, since $\phi,\psi \in \IA_N(\ZZ/3\ZZ)$, the element $g$ is fixed by $\langle \Phi,\Psi \rangle$. Note that, since $\phi=\psi$, we have $\Phi\Psi^{-1} \in K$. Since $K_0$ is cyclic, we have $\Phi=\Psi \circ \ad_{g^m}$ for some $m \in \ZZ$. Since $g$ is fixed by $\Psi$, we see that $\Phi^k=\Psi^k \circ \ad_{g^{km}}$. Since $\Phi^k=\Psi^k$, we have $\ad_{g^{km}}=\id$, so that $g=e$ and $\Phi=\Psi$. This concludes the proof when $K_0$ is cyclic.

Suppose now that $K_0$ is a nonabelian free group. Since $\Phi^k$ is central in $\langle \Phi,\Psi \rangle$, it fixes every element $g \in K_0$. Since $\phi,\psi \in \IA_N(\ZZ/3\ZZ)$, by Theorem~\ref{Theo conjugacy class fixed}, both $\phi$ and $\psi$ fix the conjugacy class of every element in $K_0$. Thus, there exist $g_{\phi}, g_{\psi} \in F_N$ such that, for every $g \in K_0$, we have $\Phi(g)=g_{\phi}gg_{\phi}^{-1}$ and $\Psi(g)=g_{\psi}gg_{\psi}^{-1}$.

We claim that $g_{\phi}$ and $g_{\psi}$ are fixed by $\Phi^k$. We prove the result for $g_{\phi}$, the arguments for $g_{\psi}$ being similar. Since the group $K$ is normal in $\langle \Phi,\Psi \rangle$, the automorphism $\Phi$ preserves $K_0$. In particular, for every $g \in K_0$, we have $=g_{\phi}gg_{\phi}^{-1} \in K_0$. Thus, for every $g \in K_0$, we have $$\Phi^k(g_{\phi}gg_{\phi}^{-1})=\Phi^k(g_{\phi})g\Phi^k(g_{\phi}^{-1})=g_{\phi}gg_{\phi}^{-1},$$ where the first equality follows from the fact that $g$ is fixed by $\Phi^k$. Thus, for every $g \in K_0$, the element $\Phi^{k}(g_{\phi})^{-1}g_{\phi}$ commutes with $g$. Since this is true for every $g \in K_0$ and since the group $K_0$ is nonabelian, we have in fact $\Phi^{k}(g_{\phi})=g_{\phi}$.

Thus, the group $K_0'=\langle K_0,g_{\phi},g_{\psi} \rangle$ is fixed elementwise by $\Phi^k$. Since $K_0$ is a nonabelian free group and since $\phi, \psi \in \IA_N(\ZZ/3\ZZ)$, by Theorem~\ref{Theo conjugacy class fixed}, the automorphisms $\Phi$ and $\Psi$ act on $K_0'$ by a global conjugation by respectively $g_{\phi}$ and $g_{\psi}$. 

Let $g \in K_0$ be such that $\Phi=\ad_g \circ \Psi$, which exists since $\phi=\psi$. We claim that $g \in \langle g_{\psi} \rangle$. Indeed, otherwise $\langle g,g_{\psi} \rangle$ is a nonabelian free group of rank $2$. But $g_{\phi}=gg_{\psi}$, so that, for every $i \in \NN^*$, the element $\Phi^i(g_{\psi})$ can be written as a nontrivial product of elements in $\langle g_{\psi} \rangle$ and elements in $\langle g \rangle$ (recall that $\Phi$ acts on $K_0'$ by a global conjugation by $g_{\phi}$). This contradicts the fact that $\Phi^k(g_{\psi})=g_{\psi}$. 

Hence we have $g \in \langle g_{\psi} \rangle$. Thus, the element $g$ is fixed by $\Psi$. Therefore, we see that $\Psi$ commutes with $\ad_g$. Since $\Psi^{-1} \circ \Phi=\ad_g$, we see that $\Psi$ commutes with $\Phi$. This concludes the proof.
\hfill\qedsymbol

\begin{lem}\label{Lem HNN extension case}
Let $\Phi \in \Aut(F_N)$ with $\phi=[\Phi] \in \IA_N(\ZZ/3\ZZ)$. Let $a,b \in F_N$ and $k \in \NN$ be such that $$\prod_{i=0}^k \Phi^i(a)=\prod_{i=0}^k \Phi^i(b).$$ Then $a=b$.
\end{lem}

\dem Let $a,b \in F_N$ be such that $\prod_{i=0}^k \Phi^i(a)=\prod_{i=0}^k \Phi^i(b)$. Let $A=\prod_{i=0}^{k-1} \Phi^i(a)$ and let $B=\prod_{i=0}^{k-1} \Phi^i(b)$. Note that, for every $c \in \{a,b\}$, we have $$\prod_{i=0}^k \Phi^i(c)=c \Phi(C),$$ so that $$a^{-1}b=\Phi(AB^{-1}). $$ Moreover, we have $A\Phi^k(a)=B\Phi^k(b)$.
Thus, $$\Phi^k(a^{-1}b)=(\Phi^k(a))^{-1}\Phi^k(b)=(\Phi^k(a))^{-1}B^{-1}A\Phi^k(a). $$ 

Thus, the outer automorphism $\phi^{k+1}$ sends the conjugacy class of $AB^{-1}$ to the conjugacy class of $B^{-1}A$. Since $AB^{-1}=BB^{-1}AB^{-1}$, we see that $\phi^{k+1}$ sends the conjugacy class of $AB^{-1}$ to itself. Thus the orbit of $[AB^{-1}]$ under iteration of $\phi$ is periodic. By Theorem~\ref{Theo conjugacy class fixed}, we see that the conjugacy class of $AB^{-1}$ is fixed by $\phi$. Thus, we have $[AB^{-1}]=[a^{-1}b]$, that is $$[a\Phi(a) \ldots \Phi^k(a)\Phi^k(b^{-1})\ldots \Phi(b^{-1})b^{-1}]=[a^{-1}b]. $$ Since $$[a\Phi(a) \ldots \Phi^k(a)\Phi^k(b^{-1})\ldots \Phi(b^{-1})b^{-1}]=[\Phi(a) \ldots \Phi^k(a)\Phi^k(b^{-1})\ldots \Phi(b^{-1})b^{-1}a],$$ there exists $h \in F_N$ with $$[\Phi(a) \ldots \Phi^k(a)\Phi^k(b^{-1})\ldots \Phi(b^{-1})]=[ha^{-1}bh^{-1}a^{-1}b]. $$ Using this argument inductively, we see that there exist $g,h_0,\ldots, h_{k} \in F_N $ with $$ga^{-1}bg^{-1}\prod_{i=0}^kh_i\Phi^i(a^{-1}b)h_i^{-1}=e. $$ Since $\Phi$ preserves the conjugacy class of $a^{-1}b$, there exist $g,g_0,\ldots, g_{k} \in F_N $ with $$ga^{-1}bg^{-1}\prod_{i=0}^kg_ia^{-1}bg_i^{-1}=e. $$ Thus some product of conjugates of $a^{-1}b$ is trivial. By for instance~\cite[Corollary~4.5]{CohenLyndon1963}, this implies that $a^{-1}b=e$, which concludes the proof.
\hfill\qedsymbol

\begin{lem}\label{Lem Sporadic case}
Let $\phi, \psi \in \IA_N(\ZZ/3\ZZ)$ and let $k \in \NN^*$ be such that $\phi^k=\psi^k$. Suppose that there exists a $\langle \phi, \psi \rangle$-invariant sporadic free factor system $\mathcal{F}$ of $F_N$ such that, for every $[A] \in \mathcal{F}$, we have $\phi|_A=\psi|_A$. Then $\phi=\psi$.
\end{lem}

\dem Since $\mathcal{F}$ is sporadic, the Bass-Serre tree $\mathcal{S}$ of $F_N$ associated with $\mathcal{F}$ is a Grushko $(F_N,\mathcal{F})$-free splitting fixed by $\Out(F_N,\mathcal{F})$, and in particular by $\langle \phi,\psi \rangle$. Since we have $\langle \phi,\psi \rangle \subseteq \IA_N(\ZZ/3\ZZ)$, by Lemma~\ref{Lem periodic splitting is fixed}, we also have $\langle \phi,\psi \rangle \subseteq K(\mathcal{S})$. We distinguish between two cases, according to the distinct cases in the definition of a sporadic free factor system.

Suppose first that $\mathcal{F}=\{[A],[B]\}$, so that $F_N=A \ast B$. For every $\theta \in K(\mathcal{S})$, let $\Theta$ be the unique representative of $\theta$ such that $\Theta(A)=A$ and $\Theta(B)=B$. By~\cite[Proposition~4.2]{levitt2005}, the map $\theta \mapsto (\Theta|_A,\Theta|_B)$ gives an isomorphism between $K(\mathcal{S})$ and $\Aut(A) \times \Aut(B)$. Note that, by assumption on $\phi$ and $\psi$, the outer automorphism class of $\Phi|_A$ (resp. $\Phi|_B)$ is the same as the one of $\Psi|_A$ (resp. $\Psi|_B)$. Since $\phi|_A,\psi|_A \in \IA(A,\ZZ/3\ZZ)$ and $\phi|_B,\psi|_B \in \IA(B,\ZZ/3\ZZ)$, by Lemma~\ref{Lem case automorphism}, we have $\Phi|_A=\Psi|_A$ and $\Phi|_B=\Psi|_B$. Thus, we have $\phi=\psi$.

Suppose now that $\mathcal{F}=\{[A]\}$ and let $g \in F_N$ be such that $F_N=A \ast \langle g \rangle$. For every $\theta \in K(\mathcal{S})$, let $\Theta$ be the unique representative of $\theta$ such that $\Theta(A)=A$ and $\Theta(g)=gg_{\theta}$ with $g_{\theta} \in A$. By~\cite[Proposition~4.2]{levitt2005}, the map $\theta \mapsto (g_{\theta},\Theta|_A)$ gives an isomorphism between $K(\mathcal{S})$ and $A \rtimes \Aut(A)$. As in the first case, by Lemma~\ref{Lem case automorphism}, we have $\Phi|_A=\Psi|_A$. 

It now suffices to prove that $g_{\phi}=g_{\psi}$. But $g_{\phi^k}=g_{\psi^k}$. Thus, we have $$g_{\phi^k}=\prod_{i=0}^{k-1} \Phi^i(g_{\phi})=g_{\psi^k}=\prod_{i=0}^{k-1} \Psi^i(g_{\psi})=\prod_{i=0}^{k-1} \Phi^i(g_{\psi}), $$ where the last equality follows from the fact that $\Phi|_A=\Psi|_A$ and from the fact that $g_{\psi} \in A$. By Lemma~\ref{Lem HNN extension case}, we have $g_{\phi}=g_{\psi}$. This concludes the proof.
\hfill\qedsymbol

\section{Proof of Theorem~\ref{Theo intro} }\label{Section Proof}

In this section, we prove Theorem~\ref{Theo intro} in the general case. First, we need a proposition regarding stabilizers of $(F_N,\mathcal{F})$-arational trees which is of independant interest. 

\begin{prop}\label{Prop nonsporadic case}
Let $H$ be a subgroup of $\IA_N(\ZZ/3\ZZ)$ and let $\mathcal{F}$ be a nonsporadic free factor system fixed by $H$. Suppose that $H$ fixes the homothety class of an $(F_N,\mathcal{F})$-arational tree $T$ and, for every $[A] \in \mathcal{F}$, that the image of the homomorphism $H \to \Out(A)$ is abelian. Then $H$ is abelian.
\end{prop}

\dem By Lemma~\ref{Lem Tits}, it suffices to show that $H$ does not contain a nonabelian free group. 

Note that the stretching factor homomorphism $\mathrm{SF}$ associated with $\Stab([T])$ has abelian image. Thus, the group $H$ contains a nonabelian free group if and only if so does $K=\Stab_{isom}(T) \cap H$.

Let $\mathcal{S}_{T,K}$ be the canonical $(F_N,\mathcal{F})$-splitting associated with $T$ and $K$ given by Proposition~\ref{Prop existence splitting arational tree}. Let $V=V_0 \amalg V_1$ be the associated partition of the vertices, let $v \in V_0$ and let $v_1,\ldots,v_k$ be representatives of the $F_N$-orbits in $V_1$. Note that, by Proposition~\ref{Prop existence splitting arational tree}~$(2)(b)$ the set $V_0$ consists in a unique $F_N$-orbit. Up to taking a finite index subgroup of $K$, we have a homomorphism $$\Lambda \colon K \to \Out(G_v) \times \prod_{i=1}^k \Out(G_{v_i}).$$ 

By Proposition~\ref{Prop existence splitting arational tree}~$(2)(a)$, for every $i \in \{1,\ldots,k\}$, there exists $[A] \in \mathcal{F}$ such that $[G_{v_i}]=[A]$. Thus, by assumption, for every $i \in \{1,\ldots,k\}$, the image of $K$ in $\Out(G_{v_i})$ is abelian. Moreover, by Proposition~\ref{Prop existence splitting arational tree}~$(2)(b)$, the image of $K$ in $\Out(G_v)$ is trivial. Hence the image of $\Lambda$ is abelian. 

Thus $H$ contains a nonabelian free subgroup if and only if so does $\ker(\Lambda)$. But $\ker(\Lambda)$ is contained in the group of twists of $\mathcal{S}_T$. Recall that, by Proposition~\ref{Prop existence splitting arational tree}~$(1)$, every edge stabilizer of $\mathcal{S}_T$ is nontrivial. Hence, by Lemma~\ref{Lem Dehn twists central}, the group of twists of $\mathcal{S}_T$ is central in $K(\mathcal{S})$. Thus, the group $\ker(\Lambda)$ is abelian. This concludes the proof.
\hfill\qedsymbol

\bigskip

\noindent{\it Proof of Theorem~\ref{Theo intro}. } Let $\phi, \psi \in \IA_N(\ZZ/3\ZZ)$ and suppose that there exists $k \in \NN^*$ such that $\phi^k=\psi^k$. We prove that $\phi=\psi$ by induction on the rank $N$ of $F_N$. When $N=1$, we have $\IA_N(\ZZ/3\ZZ)=\{\id\}$ and there is nothing to prove. 

Suppose that $N \geq 2$ and let $H=\langle \phi,\psi \rangle$. Let $\mathcal{F}=\{[A_1],\ldots,[A_{\ell}]\}$ be a maximal proper $H$-invariant free factor system. Since $\phi, \psi \in \IA_N(\ZZ/3\ZZ)$, by Theorem~\ref{Theo free factor fixed}, every element of $\mathcal{F}$ is fixed by $H$. Therefore, for every $i \in \{1,\ldots,\ell\}$, we have a natural homomorphism $H \to \Out(A_i)$ whose image is contained in $\IA_N(A_i,\ZZ/3\ZZ)$. By induction, for every $i \in \{1,\ldots,\ell\}$, we have $\phi|_{A_i}=\psi|_{A_i}$. Hence the image of the homomorphism $$\Lambda \colon H \to \prod_{i=1}^{\ell} \Out(A_i)$$ is cyclic. Thus, if $\mathcal{F}$ is sporadic, by Lemma~\ref{Lem Sporadic case}, we see that $\phi=\psi$. 

Therefore, we may suppose that $\mathcal{F}$ is a nonsporadic free factor system. We prove that the group $H$ is abelian. Since $\IA_N(\ZZ/3\ZZ)$ is torsion free by Proposition~\ref{Prop torsion free} and since $\phi^k=\psi^k$, this will conclude the proof. By Lemma~\ref{Lem Tits}, it suffices to prove that $H$ does not contain a nonabelian free group.

We claim that $\phi$ and $\psi$ are fully irreducible relative to $\mathcal{F}$. Indeed, otherwise, we may suppose that some power of $\phi$ fixes a free factor system $\mathcal{F}'$ with $\mathcal{F}<\mathcal{F}'< \{[F_N]\}$. By Theorem~\ref{Theo free factor fixed}, the element $\phi$ fixes $\mathcal{F}'$. Thus, $\phi^k=\psi^k$ fixes $\mathcal{F}'$. By Theorem~\ref{Theo free factor fixed}, we see that $\psi$ and $H$ fix $\mathcal{F}'$. This contradicts the maximality of $\mathcal{F}$.

Hence $\phi$ and $\psi$ are fully irreducible relative to $\mathcal{F}$. Thus, by Theorem~\ref{Theo loxo free factor}, they act loxodromically on the Gromov hyperbolic free factor graph $\FF(F_N,\mathcal{F})$ relative to $\mathcal{F}$. Note that both $\phi$ and $\psi$ fix exactly two points in $\partial_{\infty} \FF(F_N,\mathcal{F})$ by hyperbolicity of $\FF(F_N,\mathcal{F})$. Since $\phi^k=\psi^k$, it follows that $\phi$ and $\psi$ fix the same two points in the Gromov boundary $\partial_{\infty} \FF(F_N,\mathcal{F})$. Thus, the group $H$ has a finite orbit in $\partial_{\infty} \FF(F_N,\mathcal{F})$. 

By Proposition~\ref{prop fix arational boundary}, $H$ has a finite index subgroup $H'$ which fixes the homothety class of an arational $(F_N,\mathcal{F})$-tree. Note that $H$ contains a nonabelian free group if and only if so does $H'$.

But the group $$\Lambda(H') \subseteq \prod_{i=1}^k \Out(A_i)$$ is cyclic by hypothesis. Thus, we can apply Proposition~\ref{Prop nonsporadic case} to show that $H'$ is abelian. Hence $H$ does not contain a nonabelian free group. This concludes the proof.
\hfill\qedsymbol

\bigskip

Following Kontorovi\v{c}~\cite{Kontorovic1946} (see also~\cite{Baumslag1960}), Theorem~\ref{Theo intro} implies that the group $\IA_N(\ZZ/3\ZZ)$ is an \emph{$R$-group}. The following properties then follow for instance from~\cite[Proposition~2.2]{FayWalls1999}.

\begin{coro}\label{Coro R group}
The group $\IA_N(\ZZ/3\ZZ)$ satisfies the following properties.

\medskip

\noindent{$(1)$ } For every $\phi \in \IA_N(\ZZ/3\ZZ)$, the normalizer of $\langle \phi \rangle$ in  $\IA_N(\ZZ/3\ZZ)$ coincides with its centralizer.

\medskip

\noindent{$(2)$ } For every $\phi, \psi \in \IA_N(\ZZ/3\ZZ)$, if there exist $m,n \in \ZZ^*$ such that $\phi^m$ and $\psi^n$ commute, then $\phi$ and $\psi$ commute.
\hfill\qedsymbol
\end{coro}

\section{Abelian subgroups of $\mathrm{IA}_N(\mathbb{Z}/3\mathbb{Z})$}\label{Section abelian}

In this section, we study some properties of abelian subgroups of $\IA_N(\ZZ/3\ZZ)$. These properties are mainly consequences of Proposition~\ref{Prop nonsporadic case} and Corollary~\ref{Coro R group}.

\begin{theo}\label{Theo Normalizer equals centralizer}
Let $N \geq 2$ and let $H$ be an abelian subgroup of $\IA_N(\ZZ/3\ZZ)$. The normalizer of $H$ in $\IA_N(\ZZ/3\ZZ)$ is equal to its centralizer.
\end{theo}

\dem Let $\phi \in N_{\IA_N(\ZZ/3\ZZ)}(H)$. In order to prove Theorem~\ref{Theo Normalizer equals centralizer}, it suffices to prove that the group $H'=\langle \phi,H\rangle$ is abelian.

We prove the result by induction on $N$, the case $N=1$ being immediate. Suppose that $N \geq 2$ and let $\mathcal{F}$ be a maximal proper $H'$-invariant free factor system of $F_N$. We distinguish between two cases, according to the nature of $\mathcal{F}$.

\medskip

\noindent{\bf Case A. } Suppose that $\mathcal{F}$ is nonsporadic. 

\medskip

Since $\mathcal{F}$ is maximal, by Theorem~\ref{Theo fully irreducible contained}, the group $H'$ contains a fully irreducible outer automorphism relative to $\mathcal{F}$. By Theorem~\ref{Theo loxo free factor}, the group $H'$ contains a loxodromic element of the relative free factor graph $\FF(F_N,\mathcal{F})$, hence has unbounded orbits in $\FF(F_N,\mathcal{F})$. 

Suppose first that $H$ contains a loxodromic element $\psi$ of $\FF(F_N,\mathcal{F})$. Then $\psi$ has exactly two periodic points in $\partial_{\infty}\FF(F_N,\mathcal{F})$, which are its attracting and repelling fixed points. Since $H$ is abelian, the group $H$ has also at most two finite orbits in $\partial_{\infty}\FF(F_N,\mathcal{F})$, and such an $H$-finite orbit exists. This implies that $H'$ has a finite orbit in $\partial_{\infty}\FF(F_N,\mathcal{F})$. By Proposition~\ref{prop fix arational boundary}, the group $H'$ has a finite index normal subgroup $H_0'$ which fixes the homothety class of an $(F_N,\mathcal{F})$-arational tree $T$. 

Let $H_0$ be a characteristic subgroup of $H$ of finite index and contained in $H_0'$ (it exists since $H$ is finitely generated by Theorem~\ref{Theo virtually abelian is abelian}). Let $k \in \NN^*$ be such that $\phi^k \in H_0'$. Then the element $\phi^k$ normalizes $H_0$.

Note that, since $\langle H_0,\phi^k \rangle \subseteq \IA_N(\ZZ/3\ZZ)$, by Theorem~\ref{Theo free factor fixed}, for every $[A] \in \mathcal{F}$, we have a homomorphism $\langle H_0,\phi^k \rangle \to \Out(A)$ whose image is abelian by induction hypothesis. Hence, by Proposition~\ref{Prop nonsporadic case}, the group $\langle H_0,\phi^k \rangle$ is abelian. By Corollary~\ref{Coro R group}~$(2)$, we see that $\phi$ is in the centralizer of $H$.

Suppose now that $H$ does not contain any loxodromic element of $\FF(F_N,\mathcal{F})$. Since $H$ is a normal subgroup of $H'$ and since $H'$ has unbounded orbits, by for instance~\cite[Proposition~4.2]{HorbezWade20}, the group $H$ has a (not necessarily unique) finite orbit in $\partial_{\infty} \FF(F_N,\mathcal{F})$. 

By Proposition~\ref{prop fix arational boundary}, the group $H$ contains a finite index subgroup $H_0$ which fixes the homothety class of an $(F_N,\mathcal{F})$-arational tree $T$. Since $H$ is abelian, it is finitely generated by Theorem~\ref{Theo virtually abelian is abelian}. Hence we may suppose that $H_0$ is also a characteristic subgroup of $H$, so that $H' \subseteq N(H_0)$. 

By Lemma~\ref{Lem stratching factor cyclic}, since $H_0$ does not contain any fully irreducible element of $F_N$ relative to $\mathcal{F}$, we have in fact $H_0 \subseteq \Stab_{isom}(T)$. Let $\mathcal{S}_{T,H_0}$ be the $(F_N,\mathcal{F})$-splitting of $F_N$ given by Proposition~\ref{Prop existence splitting arational tree}. The splitting $\mathcal{S}_{T,H_0}$ is fixed by the normalizer of $H_0$ in $\Out(F_N,\mathcal{F})$. In particular, it is fixed by $H'$. Let $V=V_0 \amalg V_1$ be the associated partition of the vertices, let $v \in V_0$ and let $v_1,\ldots,v_k$ be representatives of the $F_N$-orbits in $V_1$. Note that, by Proposition~\ref{Prop existence splitting arational tree}~$(2)(b)$ the set $V_0$ consists in a unique $F_N$-orbit. 

Recall that, by Proposition~\ref{Prop existence splitting arational tree}~$~(2)(a)$, for every $i \in \{1,\ldots,k\}$, we have $[G_{v_i}] \in \mathcal{F}$. Since $H' \subseteq \IA_N(\ZZ/3\ZZ)$, Theorem~\ref{Theo free factor fixed} imples that, for every $i \in \{1,\ldots,k\}$, the group $H'$ fixes the conjugacy class of $G_{v_i}$. Since $V_0$ consists in a unique $F_N$-orbit, we see that $H' \subseteq K(\mathcal{S}_{T,H_0})$. Therefore, we have a homomorphism $$\Lambda_0 \colon H' \to \Out(G_v) \times \prod_{i=1}^k \Out(G_{v_i}). $$ By induction and Proposition~\ref{Prop existence splitting arational tree}~$(2)(a)$, for every $i \in \{1,\ldots,k\}$, the homomorphism $H' \to \Out(G_{v_i})$ has abelian image.

By Proposition~\ref{Prop existence splitting arational tree}~$(2)(b)$, the image of the homomorphism $H_0 \to \Out(G_v)$ is trivial. By~\cite[Lemma~5.5]{HorbezWade20}, this implies that the image of the homomorphism $H \to \Out(G_v)$ is also trivial. Since $H'=\langle H,\phi \rangle$, the homomorphism $H' \to \Out(G_v)$ has cyclic image. Thus $H'$ contains a nonabelian free group if and only if $\ker(\Lambda_0)$ contains a nonabelian free group.

By Proposition~\ref{Prop existence splitting arational tree}~$(1)$, every edge stabilizer of $\mathcal{S}_{T,H_0}$ is nontrivial. By Lemma~\ref{Lem Dehn twists central}, the kernel $\ker(\Lambda_0)$ is abelian. This shows that $H'$ does not contain a nonabelian free group and, by Lemma~\ref{Lem Tits}, that $H'$ is abelian. This proves Theorem~\ref{Theo Normalizer equals centralizer} when $\mathcal{F}$ is nonsporadic.

\medskip

\noindent{\bf Case B. } Suppose that $\mathcal{F}$ is sporadic. 

\medskip

The group $H'$ then fixes the splitting $\mathcal{S}$ associated with the Bass-Serre tree of $F_N$ relative to $\mathcal{F}$. Let $V$ be a set of representatives of the orbits of vertices of $S \in \mathcal{S}$. 

Since $H' \subseteq \IA_N(\ZZ/3\ZZ)$, by Lemma~\ref{Lem periodic splitting is fixed}, we have $H' \subseteq K(\mathcal{S})$. Thus, we have a homomorphism $$\Lambda \colon H' \to \prod_{v \in V}\Out(G_v)$$ whose kernel is $T(\mathcal{S}) \cap H'$. Since $\mathcal{S}$ is the Bass-Serre tree of $F_N$ associated with $\mathcal{F}$, for every $v \in V$, we have $[G_v] \in \mathcal{F}$. By induction hypothesis, the image of $\Lambda$ is abelian.

Recall that, by Lemma~\ref{Lem Tits}, in order to prove that $H'$ is abelian, it suffices to prove that $H'$ does not contain a nonabelian free group. Note that, since the image of $\Lambda$ is abelian, the group $H'$ contains a nonabelian free group if and only if $T(\mathcal{S}) \cap H'$ contains a nonabelian free group. Note also that, since $H$ is normal in $H'=\langle \phi,H \rangle$, for every $\psi \in H'$, there exist $k_{\psi} \in \ZZ$ and $h_{\psi} \in H$ such that $\psi=\phi^{k_{\psi}}h_{\psi}$. We distinguish between several cases, according to the intersection $H' \cap T(\mathcal{S})$.

\medskip

\noindent{\bf Subcase 1. } Suppose that, for every $\psi \in T(\mathcal{S}) \cap H'$, we have $k_{\psi}=0$.

\medskip

In that case, we see that $T(\mathcal{S}) \cap H' \subseteq H$. Since $H$ is abelian, this implies that $T(\mathcal{S}) \cap H'$ is abelian. Thus, $H'$ is abelian.

\medskip

\noindent{\bf Subcase 2. } Suppose that there exists $\psi \in T(\mathcal{S}) \cap H'$ with $k_{\psi} \neq 0$ and that $H \cap T(\mathcal{S})=\{1\}$. 

\medskip

In that case, $H$ and $H' \cap T(\mathcal{S})$ are normal subgroups of $H'$ with trivial intersection. Thus, the groups $H$ and $H' \cap T(\mathcal{S})$ commute. Therefore, $\psi$ commutes with every element of $H$. Since $H$ is abelian and since $h_{\psi} \in H$, this implies that $\phi^{k_{\psi}}$ commutes with every element of $H$. By Corollary~\ref{Coro R group}~$(2)$, we see that $\psi$ commutes with every element of $H$.

\medskip

\noindent{\bf Subcase 3. } Suppose that there exists $\psi \in T(\mathcal{S}) \cap H'$ with $k_{\psi} \neq 0$ and that $H \cap T(\mathcal{S})\neq \{1\}$.

\medskip

Recall that $T(\mathcal{S})$ is isomorphic to $A \times B$, where $A$ and $B$ are free (maybe cylic) groups. Since $H$ is abelian, there exist $T_A \in A$ and $T_B \in B$ such that $H \cap T(\mathcal{S}) \subseteq \langle T_A,T_B \rangle$. 

Note that $\psi$ normalizes the groups $H \cap T(\mathcal{S})$, $A$ and $B$. Since $A$ and $B$ are free groups and since $\psi \in \mathcal{T}(\mathcal{S})$, the element $\psi$ centralizes $T_A$ and $T_B$. Since $\psi \in T(\mathcal{S})$ and since $A$ and $B$ are free groups, one of the following holds: 

\noindent{$(i)$ } we have $\psi \in \langle T_A,T_B \rangle$; 

\noindent{$(ii)$ } the element $T_A$ is trivial and $\psi \in \langle A, T_B \rangle$;

\noindent{$(iii)$ } the element $T_B$ is trivial and $\psi \in \langle B,T_A \rangle$.

Suppose that $T_A$ and $T_B$ are not trivial. In that case, Assertion~$(i)$ holds for every $\psi \in H' \cap T(\mathcal{S})$.  Thus, the group $H' \cap T(\mathcal{S})$ is contained in the abelian group $\langle T_A,T_B \rangle$. Hence the group $H' \cap T(\mathcal{S})$ does not contain a nonabelian free group. By Lemma~\ref{Lem Tits}, the group $H'$ is abelian. 

Up to exchanging the roles of $A$ and $B$, suppose that Assertion~$(ii)$ holds. Then every element $\psi \in H' \cap T(\mathcal{S})$ with $k_{\psi} \neq 0$ is contained in $\langle A, T_B \rangle$. Since Assertion~$(ii)$ holds, we have $H \cap T(\mathcal{S}) \subseteq \langle T_B \rangle \subseteq B$. Since the centralizer of $T_B$ in $A \times B$ is $\langle A, T_B \rangle$, the group $H' \cap T(\mathcal{S})$ contains a nonabelian free group if and only if $H' \cap A$ contains a nonabelian free group. 

Since $H \cap T(\mathcal{S}) \subseteq B$, the groups $H' \cap A$ and $H$ are normal subgroups of $H'$ with trivial intersection. Thus, they commute with each other. The conclusion is now similar to the one in Case~2. 

As we have ruled out every case, this concludes the proof of Theorem~\ref{Theo Normalizer equals centralizer}.
\hfill\qedsymbol

\bigskip

From Theorem~\ref{Theo Normalizer equals centralizer}, we immediately deduce the following corollary, which gives obstructions to the existence of subgroups of $\IA_N(\ZZ/3\ZZ)$ isomorphic to some semidirect products of groups.

\begin{coro}\label{Coro obstruction}
Let $H$ be a subgroup of $\IA_N(\ZZ/3\ZZ)$ isomorphic to a semidirect product $A \rtimes B$, where $A$ is abelian. Then $H$ is isomorphic to the direct product $A \times B$.

In particular, if $H$ is a subgroup of $\IA_N(\ZZ/3\ZZ)$ which fits in a short exact sequence $$1 \to A \to H \to B \to 1$$ where $A$ is abelian and $B$ is a free group, then $H$ is isomorphic to $A \times B$. 
\hfill\qedsymbol
\end{coro}

Given an element $\phi \in \IA_N(\ZZ/3\ZZ)$, we denote by $\Fix(\phi)$ the set of conjugacy classes of maximal cyclic subgroups of $F_N$ fixed by $\phi$. The following result follows from the work of Feighn and Handel~\cite{FeiHan09} (see also the work of Mutanguha~\cite{Mutanguha22} regarding centralizers of \emph{atoroidal} elements of $\Out(F_n)$). Corollary~\ref{Coro R group}~$(1)$ allows us to extend this result to normalizers.

\begin{theo}\cite{FeiHan09}\label{Prop atoroidal}
Let $\phi \in \IA_N(\ZZ/3\ZZ)$ be such that $\Fix(\phi)$ is finite. The normalizer of $\langle \phi \rangle$ in $\IA_N(\ZZ/3\ZZ)$ is abelian.
\end{theo}

\dem Let $H$ be the centralizer of $\phi$ in $\IA_N(\ZZ/3\ZZ)$. By Corollary~\ref{Coro R group}~$(1)$, it suffices to prove that $H$ is abelian. The proof is by induction on $N$, the result for $N=1$ being immediate. Let $\mathcal{F}$ be a maximal proper $H$-invariant free factor system of $F_N$. As in the proof of Theorem~\ref{Theo Normalizer equals centralizer}, we distinguish between two cases, according to the nature of $\mathcal{F}$.

Suppose that $\mathcal{F}$ is sporadic and let $\mathcal{S}$ be the splitting of $F_N$ associated with the Bass-Serre tree of $F_N$ relative to $\mathcal{F}$. Let $V$ be a set of representatives of orbits of vertices in $S \in \mathcal{S}$. As in the proof of Theorem~\ref{Theo Normalizer equals centralizer}, we have a homomorphism $\Lambda \colon H \to \prod_{v \in V}\Out(G_v)$ whose image is abelian by induction hypothesis. 

We claim that the kernel of $\Lambda$ is abelian. Indeed, by Lemma~\ref{Lem Tits}, it suffices to prove that the kernel of $\Lambda$ does not contain a nonabelian free group. Recall that the kernel of $\Lambda$ is $H \cap T(\mathcal{S})$ and that the group $T(\mathcal{S})$ is isomorphic to a direct product $G_1 \times G_2$ of two free (maybe cyclic) groups. Recall also that both $G_1$ and $G_2$ are normal subgroups of $K(\mathcal{S})$. In particular, let $a \in G_1$ and $b \in G_2$ be such that $ab \in H \cap T(\mathcal{S})$. Then both $a$ and $b$ centralizes $\phi$. Thus, $H \cap T(\mathcal{S})$ contains a nonabelian free group if and only if $H \cap G_1$ or $H \cap G_2$ contains a nonabelian free group. 

For every $i \in \{1,2\}$, the map $H \cap G_i \to F_N$ which sends $\psi \in H \cap G_i$ to its twistor is injective (see for instance~\cite[Proposition~3.1]{levitt2005}). By for instance~\cite[Lemma~2.9]{HorbezWade20}, for every $i \in \{1,2\}$ and every $\psi \in H \cap G_i$, the element $\phi$ fixes the conjugacy class of the twistor associated with $\psi$. Thus, if there exists $i \in \{1,2\}$ such that $H \cap G_i$ contains a nonabelian free group, then $\Fix(\phi)$ is infinite. Therefore, since $\Fix(\phi)$ is finite, for every $i \in \{1,2\}$, the group $H \cap G_i$ does not contain a nonabelian free group. Therefore, the kernel of $\Lambda$ does not contain a nonabelian free group. By Lemma~\ref{Lem Tits}, this shows that $H$ is abelian. This concludes the proof when $\mathcal{F}$ is sporadic.

Suppose now that $\mathcal{F}$ is nonsporadic. We claim that $\phi$ is fully irreducible relative to $\mathcal{F}$. Indeed, by maximality of $\mathcal{F}$ and Theorem~\ref{Theo fully irreducible contained}, the group $H$ contains a fully irreducible outer automorphism $\psi$ relative to $\mathcal{F}$. If $\phi=\psi$, we are done. 

Otherwise, by Theorem~\ref{Theo loxo free factor}, the element $\psi$ is a loxodromic element of the relative free factor graph $\FF(F_N,\mathcal{F})$. Thus $\psi$ has exactly two fixed points in the Gromov boundary $\partial_{\infty}\FF(F_N,\mathcal{F})$. 

Since $\phi$ commutes with $\psi$, the group $\langle \phi,\psi \rangle$ has a finite orbit in $\partial_{\infty}\FF(F_N,\mathcal{F})$. By Proposition~\ref{prop fix arational boundary}, the group $\langle \phi,\psi \rangle$ has a finite index subgroup $K$ which fixes the homothety class of an $(F_N,\mathcal{F})$-arational tree $T$. Let $k \in \NN^*$ be such that $\phi^k \in K$. 

We claim that the element $\phi^k$ is not contained in the kernel of the stretching factor homomorphism $\mathrm{SF}\colon K \to \RR_+^*$. Indeed, since $\Fix(\phi)$ is finite, by Theorem~\ref{Theo conjugacy class fixed}, so is $\Fix(\phi^k)$. Therefore, there do not exist a nonabelian free group $F$ of $F_N$ and a representative $\Phi^k$ of $\phi^k$ such that $\Phi^k(F)=F$ and $\Phi^k|_F=\mathrm{id}_F$. Thus, by Corollary~\ref{Coro nonabelian fixed}, we see that $\phi^k \notin \Stab_{isom}(T)$. This proves the claim.

Thus, we have $\mathrm{SF}(\phi^k) \neq 1$. By Lemma~\ref{Lem stratching factor cyclic}, the element $\phi^k$ is fully irreducible relative to $\mathcal{F}$. Thus, $\phi$ is fully irreducible relative to $\mathcal{F}$. This proves the claim.

By Theorem~\ref{Theo loxo free factor}, the element $\phi$ fixes exactly two points in $\partial_{\infty} \FF(F_N,\mathcal{F})$. Hence $H$ has a finite orbit in $\partial_{\infty} \FF(F_N,\mathcal{F})$. By Proposition~\ref{prop fix arational boundary}, the group $H$ virtually fixes the homothety class of an $(F_N,\mathcal{F})$-arational tree. As in the proof of Theorem~\ref{Theo Normalizer equals centralizer}, we can apply Proposition~\ref{Prop nonsporadic case} to show that $H$ is virtually abelian. By Theorem~\ref{Theo virtually abelian is abelian}, the group $H$ is abelian. This concludes the proof.
\hfill\qedsymbol

\bibliographystyle{alphanum}
\bibliography{bibliographie}

\noindent \begin{tabular}{l}
Yassine Guerch \\
Univ. Lyon \\
ENS de Lyon \\
UMPA UMR 5669 \\
46 allée d'Italie \\
F-69364 Lyon cedex 07\\
{\it e-mail: yassine.guerch@ens-lyon.fr}
\end{tabular}

\end{document}